\documentclass[11pt]{article}
\usepackage{graphicx}
\usepackage{amsmath, amsfonts, amssymb}
\newcounter{theorem}

\begin{document}
\title{The Fourier Transform of the First Derivative \mbox{of the Generalized} Logistic Growth Curve}

\author{Ayse Humeyra Bilge \footnote{Faculty of Engineering and Natural Sciences, Kadir Has University, Istanbul, Turkey. Email: ayse.bilge@khas.edu.tr } \  and Yunus Ozdemir \footnote{Department of Mathematics, Anadolu University, Eskisehir, Turkey. Email: yunuso@anadolu.edu.tr } \thanks{Corresponding author}}
\maketitle

\begin{abstract}
 In this work, we obtain the Fourier transform of the first derivative of the generalized logistic growth curve in terms of Gamma functions and we discuss special cases.
 \end{abstract}



\pagestyle{myheadings}
\markboth{The Fourier Transform of the Generalized Logistic Growth}{The Fourier Transform of the Generalized Logistic Growth}

\section{Introduction}
Sigmoidal curves are monotone increasing functions $y(t)$  with horizontal asymptotes as $t\to\pm\infty$, providing mathematical models for transitions between two stable states.  In a numerical study of gelation phenomena \cite{BP2013}, we observed that the local extrema of the higher derivatives of certain sigmoidal curves tend to accumulate at a certain point.  In the search of an analytic description of this ``critical point" we have seen that the Fourier and Hilbert
transforms were the most appropriate tools.

The standard logistic growth curve is a typical example for a sigmoidal curve with an even first derivative and well known Fourier transform properties. The generalized logistic family provides good examples for  sigmoidal curves with no symmetry but the explicit expression of their  Fourier transform is not available in the literature.  The purpose of this note is to give a detailed derivation of the Fourier transform of the first derivative of the generalized logistic family. The integrals involved in the computation of this Fourier transform can be evaluated by certain computer algebra softwares, but we believe that an explicit derivation should  be found in the mathematical literature.


The standard logistic growth curve is a sigmoidal curve which is the solution of the differential equation $y'=1-y^2, \ y(0)=0$. This equation can be solved as
\begin{equation} \label{stalog}
y(t)=\tanh(t)
\end{equation}
and  its first derivative $y'(t)={\rm sech}^2(t)$ is the well known $1$-soliton solution of the Korteweg-deVries (KdV) equation.
The generalized logistic growth curve with horizontal asymptotes at $-1$ and $1$ is a sigmoidal curve given by
\begin{equation} \label{genlog}
y(t)=-1+\frac{2}{\left[1+ke^{-\beta t}\right]^{1/\nu}},
\end{equation}
where   $k>0$, $\beta>0$ and $\nu>0$.
The sigmoidal curve  (\ref{genlog}) reduces to (\ref{stalog}) for $\nu=1$, $k=1$, $\beta=2$.

 The Fourier transform $F$ of an integrable  function $f$ is defined as
\[F(\omega)=\frac{1}{\sqrt{2\pi}}\ \int_{-\infty}^{\infty}
\ f(t) \, e^{-i\omega t}\ dt,\]
for all $\omega \in \mathbb{R}$ provided that the integral exists in the sense of Cauchy principal value \cite{PapoulisFourier}.
We can recover $f$ from the inverse transform by
\[f(t)=\frac{1}{\sqrt{2\pi}}\ \int_{-\infty}^{\infty}
\ F(\omega)\, e^{i\omega t}\ d\omega,\]
for all $t \in \mathbb{R}$.
Since the sigmoidal functions, in particular the standard and generalized logistic growths, are finite as $t\to\infty$, their first derivatives are in $L^1$ and  their Fourier transform exists.
The Fourier transform of the first derivative $f(t)=y'(t)={\rm sech}^2(t)$ of the  standard logistic growth solution is obtained easily by using the integral formula
\[\int_{-\infty}^{\infty} \frac{e^{i \omega t} }{\cosh^2(t)} dt= \frac{\pi \omega}{\sinh\left(\frac{\pi \omega}{2}\right)}\]
as
\begin{equation}\label{fourier-standard}
F(\omega)=\sqrt{\frac{2}{\pi}}\frac{(\pi\omega/2)}{\sinh(\pi\omega/2)}.
\end{equation}
The first derivative of the generalized logistic growth is  a localized pulse but it is not symmetrical and the computation  of its Fourier transform is more complicated.
 In  Section 2, equation (\ref{fourier-generalized}) we obtain this Fourier transform explicitly  in terms of Gamma functions, pointing out  certain interesting relations among these and the hypergeometric functions.

We recall that the hypergeometric function $ _2F_1$ is defined by the Gauss series as
\[
_2F_1(a,b,c;z)=\sum_{n=0}^{\infty} \frac{(a)_n (b)_n}{(c)_n} \frac{z^n}{n!}
\]
on the disk $|z|<1$ (and by analytic continuation elsewhere) where $a,b,c \in \mathbb{C}$, $c\notin \mathbb{Z^{-}}\cup \{0\}$ and the symbol $(x)_n$ (also known as Pochhammer symbol) is defined by $(x)_0=1$ and $(x)_n=x(x+1)(x+n-1)$ for $1 \leq n \in \mathbb{N}$  (see \cite{Pearson} and \cite{Beukers} for details).

\section{The Fourier Transform for the Generalized Logistic Growth}
The  first derivative of the sigmoidal curve (\ref{genlog})  is
\begin{equation*} \label{asym-log-curve-derivation}
y'(t)=f(t)=\frac{2k\beta}{\nu}\left[1+ke^{-\beta t}\right]^{-\frac{1}{\nu}-1} e^{-\beta t}.
\end{equation*}
Its  Fourier transform of $f(t)$ is defined by the integral
\begin{equation}
F(\omega)=\frac{1}{\sqrt{2\pi}} \frac{2k\beta}{\nu}
\int_{-\infty}^{\infty}  e^{-i\omega t}
e^{-\beta t }\left[1+ke^{-\beta t}\right]^{-\frac{1}{\nu}-1}\, dt.\label{asym-log-curve-fourier}
\end{equation}
The definite integral
\[I(\omega):= \int_{-\infty}^{\infty}  e^{-i\omega t}
e^{-\beta t }\left[1+ke^{-\beta t}\right]^{-\frac{1}{\nu}-1}\, dt \]
can be expressed as
\[
I(\omega)=\frac{1}{\beta} \ \int_{0}^{\infty}  u^{\frac{i \omega}{\beta}}
\ \left[1+ku\right]^{-\frac{1}{\nu}-1}\ du
\]
by setting $u=e^{-\beta t}$. $I(\omega)$ can be evaluated in terms of the hypergeometric functions using the integral equality
\begin{equation}\label{bookformula}
 \int_{0}^{\infty}  x^{\lambda-1} (1+x)^\eta (1+\alpha x)^\mu \ dx = B\left(\lambda, -\eta-\mu-\lambda\right) \times \   _2F_1\left(-\mu,\lambda,-\mu-\eta;1-\alpha\right)
\end{equation}
which holds for $|{\rm arg}(\alpha)|< \pi$, $-{\rm Re}(\mu+\eta)>{\rm Re}(\lambda)>0$ (see \cite[p.317]{Grad}).

Putting $x=u$, $\lambda=1+\frac{i \omega}{\beta}$, $\alpha=k$, $\mu=-\frac{1}{\nu}-1$ and $\eta=0$ in (\ref{bookformula})
we obtain
\begin{eqnarray*}
I(\omega)= \frac{1}{\beta} \ B\left(1+\frac{i \omega}{\beta}, \frac{1}{\nu}-\frac{i \omega}{\beta}\right) \times \   _2F_1\left(\frac{1}{\nu}+1, 1+\frac{i \omega}{\beta}, \frac{1}{\nu}+1; 1-k\right) \ ,
\end{eqnarray*}
where $B$ is the well-known Beta function. It is known that
\[
B(x,y)=\frac{\Gamma(x) \Gamma(y) }{\Gamma(x+y)}
\]
provided ${\rm Re}(x)>0$ and ${\rm Re}(y)>0$, and
\[
 _2F_1\left(\frac{1}{\nu}+1, 1+\frac{i \omega}{\beta}, \frac{1}{\nu}+1;1-k\right)=k^{-1-\frac{i w}{\beta}}
\]
since $ _2F_1(b,a,b;z)=\,  _2F_1(a,b,b;z)=(1-z)^{-b}$ (see \cite[p.556]{Abrom}). Thus we have
\begin{eqnarray}\label{aysm-log-indefinite-int}
I(\omega)  = \frac{1}{\beta}  \frac{\Gamma(1+\frac{i w}{\beta}) \Gamma(\frac{1}{v}-\frac{i w}{\beta}) }{\Gamma(1+\frac{1}{\nu})} \  k^{-1-\frac{i w}{\beta}}.
\end{eqnarray}
Substituting (\ref{aysm-log-indefinite-int}) in (\ref{asym-log-curve-fourier})
\begin{eqnarray*}
F(\omega)&=&\frac{1}{\sqrt{2\pi}} \frac{2k\beta}{\nu} \ I(\omega)\\
&=&\sqrt{\frac{2}{\pi}} \frac{k^{-\frac{i w}{\beta}}}{\nu} \frac{\Gamma(1+\frac{i w}{\beta}) \Gamma(\frac{1}{v}-\frac{i w}{\beta}) }{\Gamma(1+\frac{1}{\nu})}.
\end{eqnarray*}
and using the equality $\Gamma(1+\frac{1}{\nu})=\frac{1}{\nu} \, \Gamma(\frac{1}{\nu})$, we can express the Fourier transform of the derivative of the generalized logistic curve as
\begin{equation}\label{fourier-generalized}
F(\omega)=\sqrt{\frac{2}{\pi}} \, k^{-\frac{i \omega}{\beta}} \, \frac{\Gamma(1+\frac{i \omega}{\beta}) \Gamma(\frac{1}{\nu}-\frac{i w}{\beta}) }{\Gamma(\frac{1}{\nu})} .
\end{equation}

\section{Special Cases}
We rewrite the Fourier transform pair for the first derivative as, displaying the dependence on the parameters $k$ and $\nu$ as
\begin{eqnarray}
f(t,k,\nu)&=&\frac{2k\beta}{\nu}\left[1+ke^{-\beta t}\right]^{-\frac{1}{\nu}-1} e^{-\beta t},\cr
F(\omega,k,\nu)&=&\sqrt{\frac{2}{\pi}}\frac{1}{\Gamma(\frac{1}{\nu})} \, e^{-i\left(\frac{\ln k}{\beta}\right)\omega}
\,
 \Gamma\left(1+\frac{i \omega}{\beta}\right) \Gamma\left(\frac{1}{\nu}-\frac{i w}{\beta}\right).
 \end{eqnarray}
Substituting  $\nu=1$, $k=1$, $\beta=2$ in (\ref{fourier-generalized}) and using the property
\begin{equation*}
\Gamma(x)\Gamma(1-x)=\frac{\pi}{\sin(\pi x)}
\end{equation*}
together with $\sin(ix)=i\sinh(x)$, we get the Fourier transform of the standard logistic growth function as given in (\ref{fourier-standard}).

Differentiating $f(t,k,\nu)$ with respect to $t$ and setting it equal to zero we obtain the location of the maximum of $f(t,k,\nu)$, that we denote by $t_m(k,\nu)$ as
$$t_m(k,\nu)=\frac{\ln(k/\nu)}{\beta}.$$
We recall that a shift of the origin in the time domain by an amount $\alpha$ corresponds to the multiplication of the Fourier transform  by a factor $e^{-i\alpha \omega}$, i.e, if $F(\omega)$ is the Fourier transform of $f(t)$, then $e^{-i\alpha \omega}F(\omega)$ is the Fourier transform of $f(t-\alpha)$.   Thus the parameter $k$ has the effect of shifting the origin of the time axis.
For $k=1$, and $\nu=1$, the peak of the first derivative of the standard logistic growth is located at $t=0$.  For $k=1$, and $\nu>1$, the peak shifts left while for  $\nu<1$ it shifts right.

If $\nu=1/n$, where $n$ is a positive integer greater than $1$, we use the property $\Gamma(1+x)=x\,\Gamma(x)$ to express
$F(\omega,1,\frac{1}{n})$ in terms of $F(\omega,1,1)$ as
\begin{equation*}
F(\omega,1,1/n)=\frac{1}{\Gamma(n)}(1-i\omega)(2-i\omega)\cdots (n-1-i\omega)  \ F(\omega,1,1).
\end{equation*}
This expression is a polynomial multiple of the standard logistic growth.  Since the Fourier transform of the $n$th derivative of $f(t)$ is $(i\omega)^nF(\omega)$, it follows that $f(t,1,1/n)$ is a polynomial in the derivatives of $f(t)$.

\begin{figure}[h!]
\centering\begin{minipage}{.48\textwidth}
\centering\includegraphics[width=0.99\textwidth]{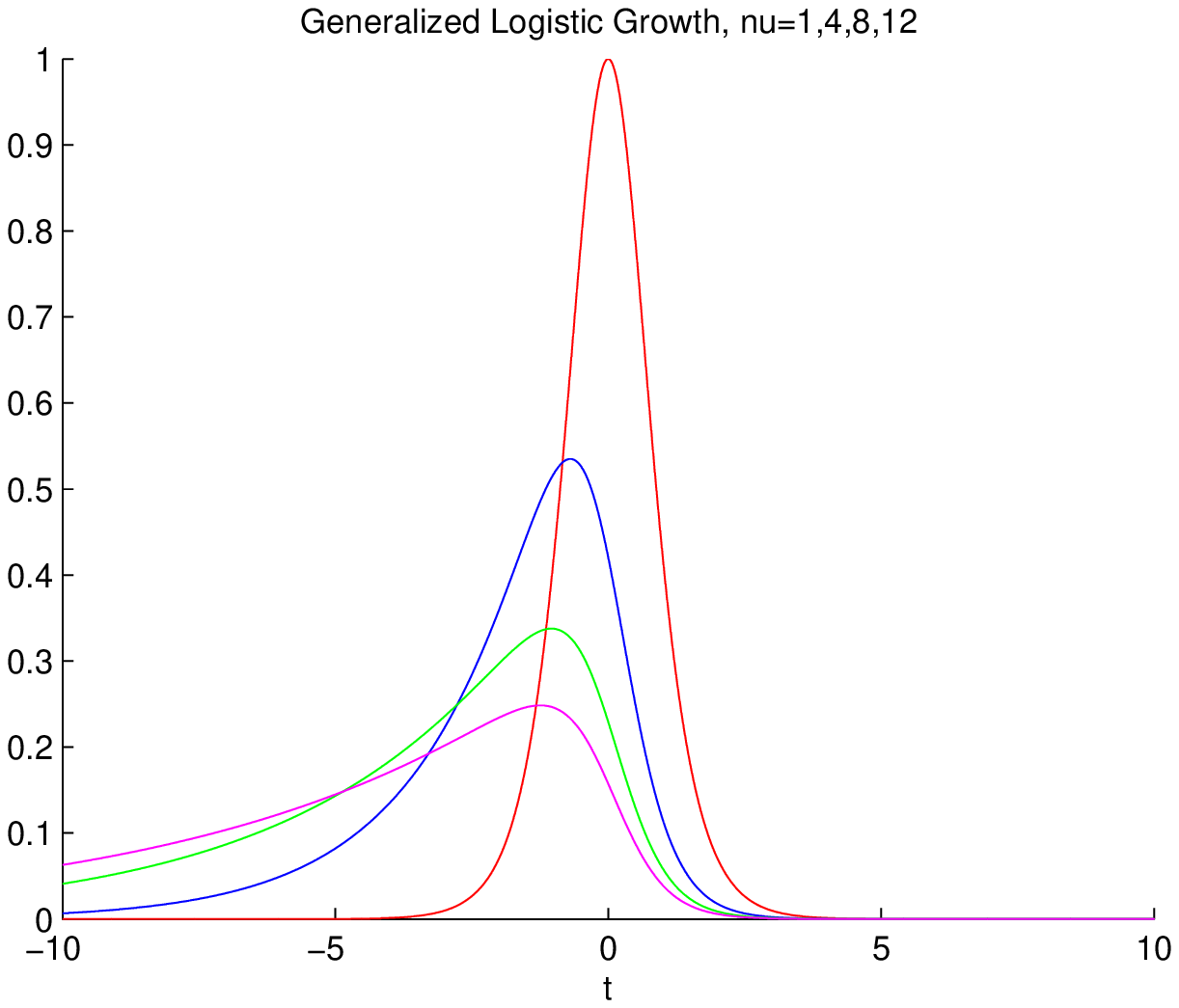}
\end{minipage}%
\centering\begin{minipage}{.48\textwidth}
\centering\includegraphics[width=0.99\textwidth]{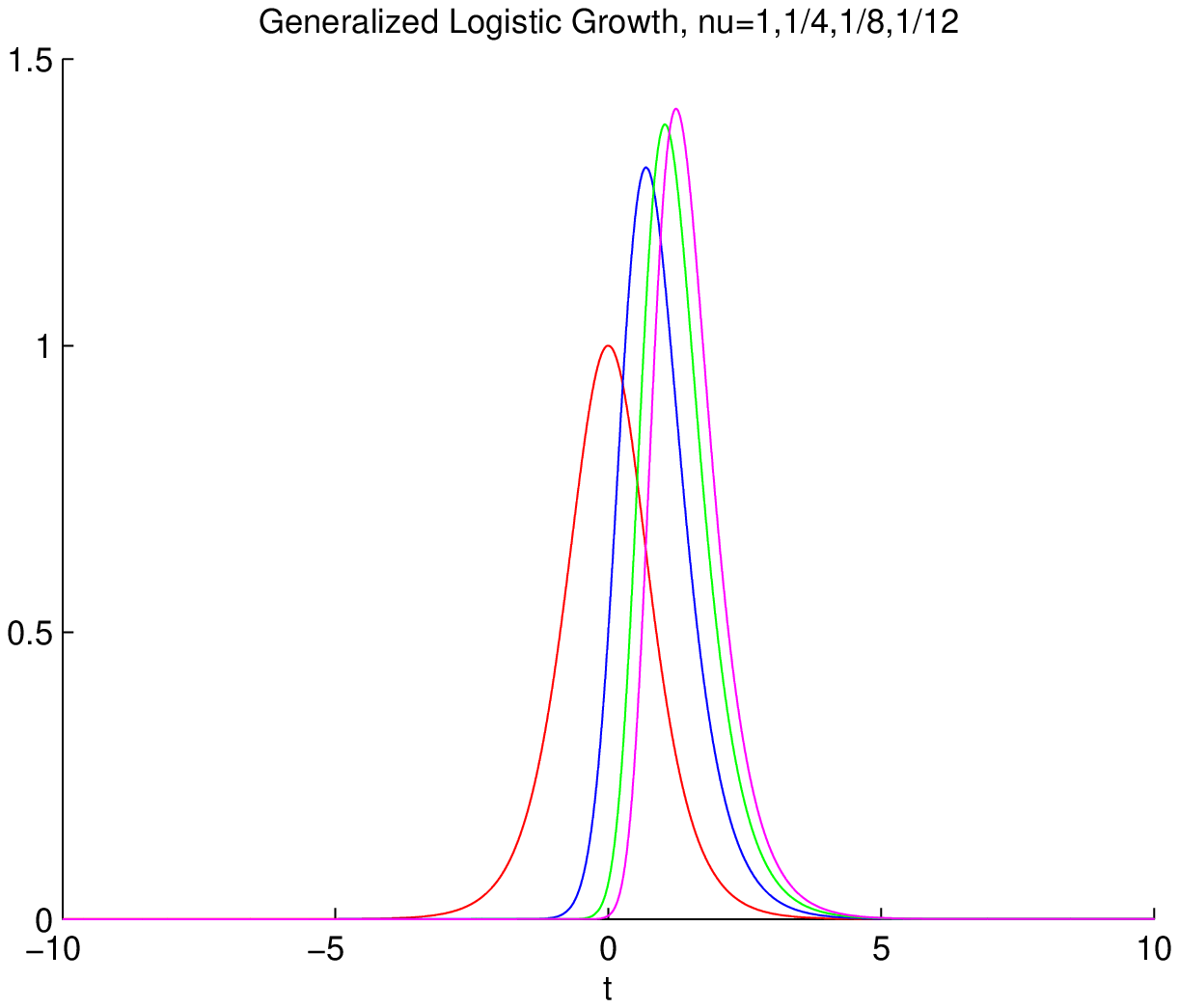}
\end{minipage}\\
   \caption{\baselineskip 12pt Time domain plots of the first derivative of the generalized logistic growth: a) for $1/\nu=1,1/4,1/8,1/12$ (left).\, b) for $1/\nu=1,4,8,12$ (right).
}\label{figtimedomain}
\end{figure}

\begin{figure}[h!]
\centering\begin{minipage}{.48\textwidth}
\centering\includegraphics[width=0.99\textwidth]{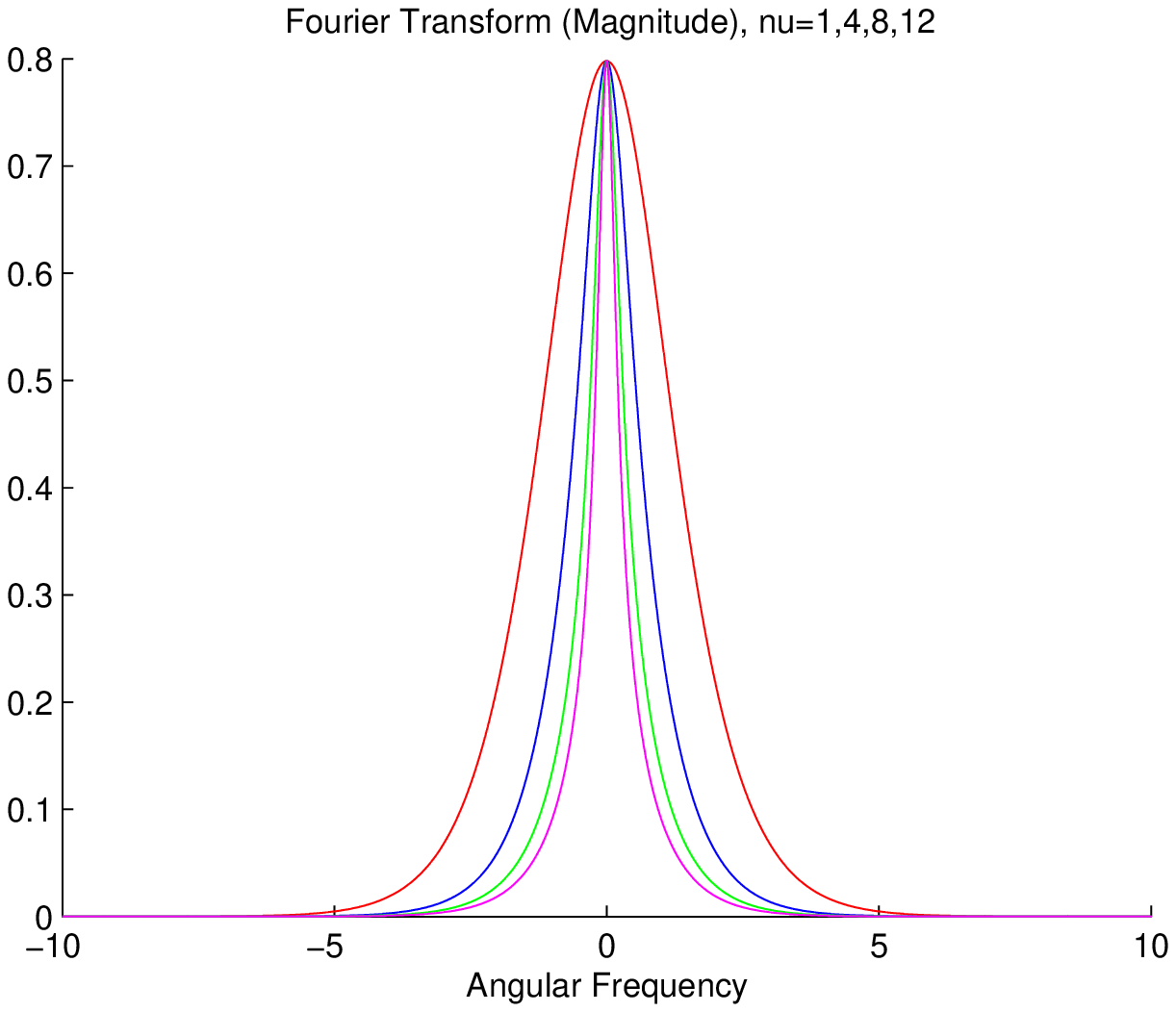}
\end{minipage}%
\centering\begin{minipage}{.48\textwidth}
\centering\includegraphics[width=0.99\textwidth]{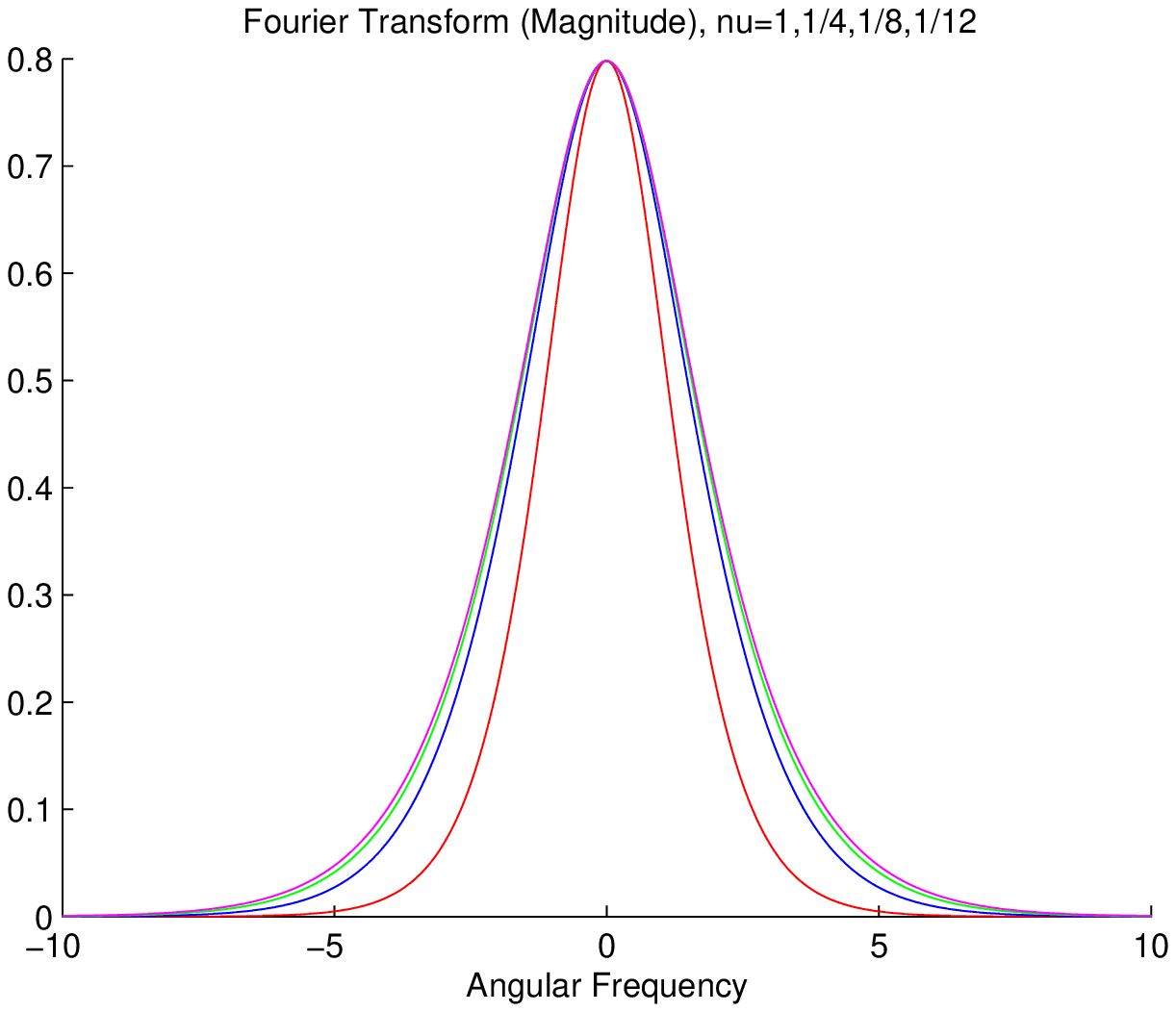}
\end{minipage}\\
   \caption{\baselineskip 12pt The magnitude of the Fourier transform: a) for $1/\nu=1,1/4,1/8,1/12$ (left). \,
   b) for $1/\nu=1,4,8,12$ (right).
}\label{figmagnitudes}
\end{figure}

\begin{figure}[h!]
\centering\includegraphics[width=0.7\textwidth]{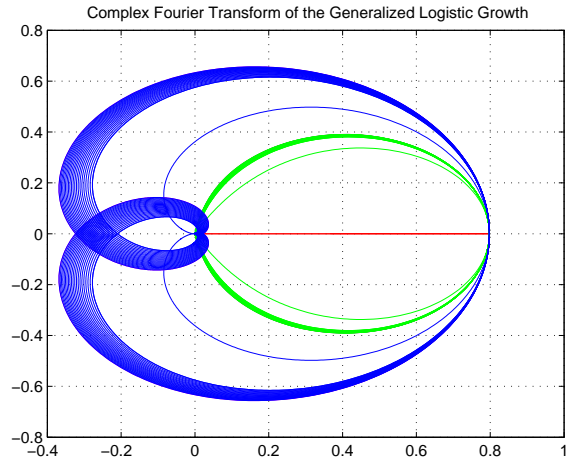}
   \caption{Parametric plot of the complex Fourier transform of the first derivative of the generalized logistic growth (from inside to out) for $\nu=1$, $\nu=4$, $\nu$ from $8$ to $12$, $\nu=1/4$ and $\nu$ from $1/8$ to $1/12$ respectively.}
\label{figparametricplot}
\end{figure}

For arbitrary values of  $\nu$, the (complex) Gamma function with complex arguments can be computed numerically.  In our case, as we are interested in the Fourier transform $F(\omega)$ for fixed values of the parameters,  we need to obtain the graphs of the real and imaginary parts of $F(\omega,1,\nu)$ on vertical lines in the complex plane.
It is known that
the Gamma functions falls of faster than any polynomial in the imaginary direction, it follows that the Fourier transform of all higher derivatives are rapidly decreasing functions.
For the cases $\nu=1/n$ and $\nu=n$, we present the plot of the first derivatives of the generalized logistic growth in Figure~\ref{figtimedomain} and the plot of the magnitudes of the Fourier transform of the first derivatives in Figure~\ref{figmagnitudes} for the values of $n=1,4,8,12$ (with $k=1, \beta=2$).
 By continuity of the Gamma function with respect to real part of its argument, the parametric plot of the complex Fourier transform for  $n<1/\nu<n+1$ fill the region between the curves corresponding to the integer values of $1/\nu$, as shown in Figure~\ref{figparametricplot}.

\clearpage
\bibliographystyle{amsplain}

\end{document}